\documentclass[11pt]{amsart}

\usepackage[dvips]{epsfig}
\usepackage{bbm,amssymb,mathabx}

\newcommand{\eps}{\varepsilon}

\def\calS{\mathcal{S}}

\def\sign{{\rm sign}}

\def\ra{\rangle}

\def\eps{{\varepsilon}}
\def\les{\lesssim}

\def\nn{\nonumber}
\def\calH{{\mathcal H}}

\def\Ran{\mathrm{Ran}}
\def\mc{\mathcal}
\def\Z{{\mathbb Z}}

\def\R{\mathbb{R}}

\def\mid{\:|\:}

\def\sign{\mathrm{sign}}

\def\f{\frac}

\def\calF{\mathcal{F}}
\def\mes{\mathrm{mes}}

\newcommand{\bmat}{\begin{pmatrix}}
\newcommand{\emat}{\end{pmatrix}}

\def\Sph{\mathbb{S}}

\def\dd{\,d}
\newcommand{\one}{\mathbbm{1}}
\def\mes{\mathcal{M}}
\def\embed{\hookrightarrow}
\newcommand{\W}{{\mathcal W}}
\def\norm[#1][#2]{\|#1\|_{#2}}
\newcommand{\B}{{\mathcal B}}
\newcommand{\K}{\Kato}
\newcommand{\Kato}{{\mathcal K}}
\newcommand{\oast}{\circledast}

\newtheorem{theorem}{Theorem}
\newtheorem{prop}[theorem]{Proposition}

\def\ov{\overline}
\def\restrict{\upharpoonright}

\def\half{\f12}
\def\Laplace{\Delta}
\def\Range{\mathrm{Ran}}
\newcommand{\EQ}[1]{\begin{equation}\begin{split} #1 \end{split}\end{equation}}
\def\trip{|\!|\!|}
\def\what{\widehat}

\begin{document}

\title[Wave operators, restriction, Wiener theorems]{Intertwining wave operators, Fourier restriction, and Wiener theorems} 
\author{
W.\ Schlag 
}

\address{ The University of Chicago, Department of Mathematics, 5734 South University Avenue, Chicago, IL 60637, U.S.A.}

\thanks{The author thanks the organizers of the Kato Centennial Conference in September of 2017 for their kind invitation, and the University of Kyoto, RIMS at Kyoto, and the University of Tokyo, for their support. The author also thanks the Institute for Advanced Study, Princeton, for its hospitality during the 2017-18 academic year.  The author was partially supported by the NSF, DMS-1500696. He is grateful to Burak Erdo\u{g}an,  Rui Han,  Marius Lemm, and Kenji Yajima for comments on an earlier version of this paper. }

\maketitle

\section{ Spectral and Scattering Theory   }

\subsection{Introduction}

This paper is an expanded version of the author's talk at the Tosio Kato centennial which took place in Tokyo, Japan, in the summer of 2017. Tosio Kato's contributions  to operator theory in general, and the spectral theory of Schr\"{o}dinger operators in particular, are monumental and we cannot do justice to them in this brief survey article. The purpose here is rather to highlight certain developments which build upon his work and would not have been possible without it. We will confine ourselves strictly to two body Schr\"{o}dinger equations. This being said, this survey is by no means exhaustive even within that more narrowly defined scope. The choices of topics is limited to the problem of asymptotic completeness, its relation to the Fourier restriction theory (albeit in the simpler Stein-Tomas incarnation which does not rely on deep geometric considerations of the Kakeya type), and finally, the $L^{p}$ theory of the wave operators initiated by Kenji Yajima\footnote{To mention some omissions:  Kato's work on trace class perturbations~\cite{Katobook,reesim}, Mourre theory~\cite{mou, HSS}, the Enss method~\cite{ens,DG}, microlocal techniques in scattering theory~\cite{ike1,II1,II2,ik1,ik2,DG}.}. The author hopes to present  these topics from a fairly general point of view, with the goal of pointing out the relevance of Fourier restriction phenomena to the study of asymptotic completeness and wave operators. The most recent results in this survey establish a structure formula for the intertwining wave operators in $\R^{3}$, obtained jointly with Marius Beceanu in 2016. 

\smallskip 

On the one hand, these results also rely on Stein-Tomas type Fourier restriction results since they depend on the more recent form of the classical Agmon-Kato-Kuroda theory due to M.~Goldberg, A.~Ionescu, and the author, \cite{golsch}, \cite{ionschlag}. On the other hand, they  serve to illustrate the power of  Wiener-type inversion theorems in Banach algebras, as applied to spectral theory.  This device was introduced into  spectral and scattering theory in Beceanu's 2009 Ph.D.\ thesis \cite{bec}. It is a   powerful device which allows one to sum otherwise divergent Born series expansions. The crucial  non-vanishing condition for these Wiener theorems is provided by some form of the limiting absorption principle, in other words, the invertibility of the Birman-Schwinger operator for positive energies (whereas for zero energy it is guaranteed via the assumption that there are no zero energy eigenvalues or resonance\footnote{The latter refers to a nontrivial solution $\psi$ of $H\psi=0$ which is not in $L^{2}$, but satisfies other types of dimension-dependent boundedness conditions,  assuming suitable  decay of the potential $V$. In dimension $d=1$ a resonance function is required to be bounded -- hence the free Laplacian exhibits a $0$ energy resonance --  and in dimension $d=3$ a resonance function decays at the rate $|x|^{-1}$. In dimension $d=2$,  there are two different kinds of resonance functions, namely $s$-waves and $p$-waves, see~\cite{EG, EGG}.  If $d>4$ a resonance at $0$ energy does not occur since the Newton potential is $L^{2}$ at $\infty$. Independently of the dimension this  $0$ energy obstruction is characterized via the Laurent expansion of the resolvent near $z=0$, cf.~\eqref{eq:symmRI} and the discussion following it.}).  This invertibility condition already appears in the classical Agmon-Kato-Kuroda theory from the 1960s and 70s. However, the form in which it appears here falls outside the scope of this older theory, which is based on weighted $L^{2}$ spaces. The form in which it arises in the aforementioned structure theorems depends on $L^{p}$ spaces rather than weighted $L^{2}$,  leading directly into Fourier restriction techniques.  The most delicate aspect of using Wiener theorems in spectral theory is to find the right spaces and algebras. We will give some indication of this in Section~\ref{sec:wiener}. 

\subsection{Wave operators and asymptotic completeness} 

 Let {   $V$} be a real-valued potential in {   $\R^d$}, bounded, and sufficiently decaying,  and set {   $H:=-\Delta+V$, $H_0:=-\Delta$}.  Define the wave operators 
 \EQ{\label{eq:W}
 W_{\pm}:=\lim_{t\to\mp\infty} e^{itH} e^{-itH_0}
 } 
These limits are known to exist  in the {    strong} {   $L^2$}-sense, provided $V$ has sufficient decay. To illustrate this, suppose    $d\ge 3$, $f\in L^1\cap L^2(\R^d)$, and assume the potential $V$ lies in~$L^2$. Then  we have 
 \begin{align}
 W_{\pm} f &=  f \mp i \int_0^\infty e^{itH} V e^{-itH_0} f \, dt \label{eq:cook} \\
  \int_1^\infty    \big \|e^{itH} V e^{-itH_0} f \big\|_2 \, dt   & \le \int_1^\infty   \| V\|_2 \| e^{-itH_0} f  \|_\infty \, dt \nn  \\
  &\les \| V\|_2  \int_1^\infty t^{-\frac{d}{2}} \|f\|_1\, dt  < \infty   \nn  
 \end{align} 
 using the pointwise decay of the free Schr\"{o}dinger evolution. Thus, the integral in the first line converges absolutely in the~$L^{2}$ norm (this is called {\em Cook's method}). 
By  unitarity of the Schr\"{o}dinger evolution, and the density of {   $L^1\cap L^2(\R^d)$} in {   $L^2$}, we conclude that  the limit exists for all {   $f\in L^2$} and that
{   $W_{\pm}$} are isometries.  Note that the condition $V\in L^{2}$ is in general not optimal in terms of decay at infinity for the existence of wave operators. 

It is clear that $e^{isH} W_{\pm} =  W_{\pm} e^{isH_{0}}$ for all $s\in\R$, and therefore by the Fourier transform also 
   $$f(H) W_{\pm} = W_{\pm} f(H_0)$$ for Schwartz functions~$f$. 
   This is precisely the {\em intertwining property}. 
   By general properties of isometries  we conclude that 
   \EQ{\label{eq:intertwine}
   f(H)  P = f(H) W_{\pm}W_{\pm}^* = W_{\pm} f(H_0)W_{\pm}^*,
   }
where {   $P$} is the  orthogonal projection onto {   $\Range(W_\pm)$}.  By the dispersive decay of the free Schr\"{o}dinger evolution, one further has  {     $\Ran(W_{\pm})\perp L^2_{pp}$}. The latter is the subspace spanned by the eigenfunctions of {   $H$}.  In fact, the more precise inclusion {     $\Ran(W_{\pm})\subset L^2_{ac}(\R^d)$}  (the absolutely continuous subspace) holds. This is implied by~\eqref{eq:intertwine}, which in turn only depends on the existence of the wave operators as  strong limits. 
 
 The fundamental 
 {\em  Asymptotic Completeness}  property goes beyond this and  states that      $\Ran(W_{\pm}) = L^2_{ac}(\R^d)$ and $L^2_{sc}=\{0\}$ (the singular continuous subspace). 
 The {    Agmon-Kato-Kuroda theory of the 1960s, and early 70s established that  the  {\em    short range condition}   
 \EQ{\label{eq:short}
   |V(x)|\les \langle x\ra^{-1-\eps}
   }
    guarantees this property, and by an earlier theorem of Kato there are {    no embedded eigenvalues} in the continuous spectrum {   $[0,\infty)$} for such potentials, cf.~ the classical papers \cite{kato, kaku, agmon, kur1, kur2, kur3} and the books~\cite{reesim,yaf1,yaf2,eskin}. 

 This theory  is based on the {\em Trace Lemma}: for   $\gamma>\frac12$} 
 \EQ{\label{eq:tracelem}
 \| \hat f \restrict S\|_{L^2(S)} \le C(\gamma,S) \| \langle x\ra^{\gamma} f\|_{L^2(\R^d)},
 } for all Schwartz functions~$f$,  where {   $S\subset \R^d$} is a smooth compact hyper-surface. 
 For the proof one  straightens the surface locally into a plane,  and applies the estimate $\|\hat{f}\|_{\infty}\le \|f\|_{1}$ and Cauchy-Schwarz.   Define the {    restriction operator}
  {   $ \rho f:=  \hat{f}\restrict S$}.  Then {   $ \rho^* g = \widehat{g\sigma_S}$, $\rho^*\rho\; f = \widehat{\sigma_S}\ast f$}.     The aforementioned trace lemma is therefore equivalent with the 
  following weighted {   $L^2$} bound:  
  \EQ{\label{eq:dualTL}
  \|  w\rho^*\rho \, w\,  f  \|_2 \le C(\eps,S) \|f\|_2,
  } 
  where $w(x)= \langle x\ra^{-\frac12-\eps}$. 
  
  \subsection{ Limiting Absorption Principle }

The bound \eqref{eq:dualTL} holds 
 for the imaginary parts of the free resolvents since 
 {     
 \[
 [(-\Delta - (\lambda^2 +  i0))^{-1} - (-\Delta - (\lambda^2 -  i0))^{-1}]f =  c\lambda^{-1} \widehat{\sigma_{\lambda \Sph^{d-1}}}\ast f 
 \]
 }
 The {\em      Limiting Absorption Principle} states that \eqref{eq:dualTL} 
 remains valid for the full resolvent, viz.
 \EQ{\label{eq:LAP}
  \| w (-\Delta - (\lambda^2 +  i0))^{-1} w \, f \|_2 \le C(\lambda)  \|f\|_2,
 }
 where $C(\lambda)\to0$ as $\lambda\to\infty$. 
 
 The first step towards establishing \eqref{eq:LAP} for $H$ rather than~$H_{0}$ is the 
 {\em resolvent identity}:    
\EQ{\label{eq:RI}
R(\lambda) &= (H-(\lambda^2 +  i0))^{-1} = R_0(\lambda)  - R_0(\lambda)  V R(\lambda)  = \\
&= ... = R_0(\lambda)  - R_0(\lambda) VR_0(\lambda)  + R_0(\lambda) VR_0(\lambda) VR_0(\lambda)  - \ldots 
 } 
 If {$V$}  is short range and    small, then by means of this expansions {   $R(\lambda)$}  inherits the limiting absorption principle. Indeed, split {   $V=|V|^{\frac12} \, \sign(V)|V|^{\frac12}= |V|^{\frac12} \, U$} and note that $|V|^{\frac12}$ has the decay required by~$w$ above. Therefore, the infinite series is summable in the corresponding weighted $L^{2}$ norm. 
 
If $V$ is large, then one cannot sum the infinite series. Instead, we treat the first line of~\eqref{eq:RI} as an implicit equation for the resolvent~$R(\lambda)$, or equivalently express the resolvent in the symmetric form 
  \EQ{\label{eq:symmRI}
  R(\lambda) = R_0(\lambda)- R_0(\lambda)|V|^{\frac12} (I + UR_0(\lambda)|V|^{\frac12}  )^{-1} U R_0(\lambda). 
  }
  The main inversion problem to be solved now is that of the Birman-Schwinger operator $I+UR_0(\lambda)|V|^{\frac12}$.  The main conclusion of Agmon-Kato-Kuroda theory is that this operator does have an inverse on $L^{2}$ for all $\lambda>0$. The argument proceeds via compactness of $UR_0(\lambda)|V|^{\frac12}$, the Fredholm alternative, and the realization that the obstruction to invertibility lies with embedded eigenvalues of $H$ (which do not exist). The characterization of obstructions is the most delicate step in the argument and requires showing that embedded resonances are necessarily eigenvalues. 
  
  Zero energy $\lambda=0$ is special and $I+UR_0(0)|V|^{\frac12} $ may be invertible on $L^{2}$ or not. The latter case is equivalent to zero energy being an eigenvalue or a resonance. In the context of classical spectral theory such as asymptotic completeness and Fourier expansions via generalized eigenfunctions with an associated Plancherel theorem, the issue of zero energy eigenvalue or resonance is irrelevant. Loosely speaking, this means that a zero energy obstruction does not affect that $L^{2}$ theory. However, for questions pertaining to $L^{p}$ with $p\ne 2$ (such as dispersive decay of the Schr\"{o}dinger evolution of~$H$ or Yajima's $L^{p}$ theory of the wave operators) the behavior of zero energy  has a profound effect as we will see below. 

While the Fredholm approach to the limiting absorption principle is indirect and thus noneffective, alternatives exist which allow for quantitative control of the constants, see~\cite{RT}. 

\section{Fourier restriction}

\subsection{Stein-Tomas theorem}

In contrast to the trace lemma~\eqref{eq:tracelem} which does not take the curvature of the hyper-surface into account, one has then following classical 
{ \em    Stein-Tomas} theorem: 

\begin{theorem} If {   $S$} has nonzero Gaussian curvature, then 
    \EQ{\label{eq:SteinTomas}
    \| \hat f \restrict S\|_{L^2(S)} \le C \|  f\|_{L^{p_{d}}(\R^d)}, \qquad p_{d}={(2d+2)}/{(d+3)}.
    }
 \end{theorem}
 
 To motive this result, note the 
 trivial bound: {   $$\|\hat f \restrict S\|_{L^2(S)} \le C \|  f\|_{L^{p}(\R^d)}, \quad p=1 $$} for any compact surface {   $S$}.  This is false if {   $p=2$} by the Plancherel theorem.
 It is natural to ask:  could there exist some {   $1<p<2$} for which this remains true? If {   $S$} (a piece of) a plane, then the answer is clearly {\bf ``no''}, since  this reduces to one variable for which we need that   $\hat{f}$ to be continuous. 
 It turns out, however,  that for {    nonvanishing curvature} the answer is { \bf ``yes''}. To see this, define the restriction operator {   $\rho f:= \hat{f}\restrict S$}. 
 Its dual is given by the inverse Fourier transform 
 $$\rho^* g=  \widecheck{ g \,\sigma_S },$$ and 
 {   $\rho^* \rho f = \widecheck{\sigma_S} \ast f$}.  The Stein-Tomas theorem is equivalent to the following bound (``factoring through $L^{2}$'')
 \EQ{\label{eq:T*T}
 T:=\rho^* \rho : L^{p_d}(\R^d) \to L^{p_d'}(\R^d)
 }
For the sake of completeness we recall the main elements of the proof. 
First, we show how to cover the range {   $p\to p'$} with {   $p<p_d$}. Write 
\[
T=\sum_j T_j, \quad T_j f= \widecheck{\sigma_S}\chi_{[|x|\simeq 2^j]} \ast f, \;\; j\ge0.
\] 
Then  
 \EQ{\label{eq:tomas}
 \| T_j\|_{1\to\infty} \les 2^{-j\frac{d-1}{2}}, \quad \| T_j\|_{2\to 2} \les 2^{jd} 2^{-j(d-1)} = 2^j
 }
 The first bound uses the decay of the Fourier transform of the surface measure. By means of stationary phase, this is  a consequence of the non-degeneracy of the second fundamental form, i.e., the non-vanishing of the Gaussian curvature: $| \widecheck{\sigma_S}(\xi)|\le C\langle \xi\ra^{-\frac{d-1}{2}}$.  By Plancherel's theorem the second bound in~\eqref{eq:tomas} reduces to the size  of the intersection of a small ball with a hyper-surface and does not use curvature. 
 
 By interpolation  {   $\|T_j\|_{p\to p'} \les 1$} where {   $\theta - (1-\theta)(d-1)/2=0$}, and {   $1/p=\theta/2+1-\theta$}. This gives exactly {   $p = (2d+2)/(d+3)=p_{d}$}. For $p<p_{d}$ one gains a convergent geometric factor $2^{-j \delta}$ with some $\delta=\delta(p)>0$. To compensate for the divergence at the critical value $p=p_{d}$, one can invoke {\em     Stein complex interpolation}, which is a method {\em for summing divergent series}.  Loosely speaking, the idea is to sum first {\em with complex weights} and then interpolate, rather than first interpolate and then sum. 
 More strictly speaking, one embeds the operator  {   $T$} into a family depending analytically on a complex parameter. 

The Stein-Tomas theorem {    is sharp}, as can be seen by the {\em Knapp example}: let {   $f=\chi_K$} be smoothed out indicator function of the cap   {   $K\subset \Sph^{d-1}$ }  of diameter~$\delta$. Then {   $|\widehat{f\, \sigma_{\Sph^{d-1}}}|$} behaves (up to tails) like  an indicator  function of a cylinder of dimensions  {   $R^2\times R\times\cdots \times R\times R$} of height {   $R^{-(d-1)}$, $R=\delta^{-1}$}.  This {     exactly balances} the inequality  {   $\| \widehat{f\,\sigma}\|_{p_d'}\les \|f\|_{L^2(\Sph^{d-1})}$}. 

\subsection{ Strichartz estimates  }

Before discussing applications of the Stein-Tomas theorem to spectral theory and the intertwining operators~\eqref{eq:W} we point out the close connection between the Fourier restriction theory and another of Tosio Kato's main interests, namely nonlinear dispersive evolution equations. To be specific, 
consider the Schr\"odinger flow {   
\EQ{ \nn
e^{-it\Delta}f(x) &= \int_{\R^d} e^{ix\cdot\xi} e^{it|\xi|^2}\, \hat{f}(\xi)\, d\xi \\
&= \int_{\R^{d+1}} e^{i(x\cdot\xi + t\tau)} \, \delta(\tau-|\xi|^2) \hat{f}(\xi) \, d\xi d\tau = (\hat{f} \,\mu )^{\check{\ }}
}
}
where {   $\mu$} is the {    measure on the paraboloid} {   $\tau=|\xi|^2$} given by {   $d\xi$}. 

By the Stein-Tomas, with increased dimension {   $d\to d+1$},  we obtain 
\EQ{ \label{eq:Lq}
\big\| e^{-it\Delta}f \big\|_{L^q_{t,x}(\R^{d+1})}  &\les  \| \hat{f} \|_{L^2(\mu)} = \|f\|_{L^2(\R^d)}, \quad q= 2+\frac{4}{d}
}
Notice an essential difference between~\eqref{eq:Lq} and the Stein-Tomas theorem. While in the latter the surface is compact, here it is not. Therefore one needs to {     scale} a compact piece of the paraboloid to the full one.  This is an example of {     many Strichartz estimates}, and similar ones hold for the wave equation, cf.~\cite{keetao}. For the latter the characteristic surface  is a {     cone} with one vanishing principal curvature, so there is a ``loss'' of one dimension. In addition, there is a singularity at the origin, which brings in Littlewood-Paley theory in order to sum the contributions coming from dyadic pieces of the cone, see~\cite{Str} for the original reference. 

To illustrate the usefulness of this type of estimate consider the 
example of an  {   $L^{2}(\R^{d})$} critical nonlinear Schr\"{o}dinger equation 
{   
\[
i\partial_{t}\psi - \Delta \psi = \pm |\psi|^{\frac{4}{d}}\psi,\quad \psi(0)=\psi_{0}\in L^{2}
\]
}
It is invariant under the scaling {   $\psi_{0}(x)\to \lambda^{\frac{d}{2}}\psi_{0}(\lambda x)$, $\psi(t,x)\to \lambda^{\frac{d}{2}}\psi(\lambda^{2}t,\lambda x)$}. 
This PDE reduces to an integral equation via Duhamel's formula, to wit
\EQ{\label{eq:duh}
\psi(t) &= e^{-it\Delta} \psi_{0} \mp i \int_{0}^{t} e^{-i(t-s)\Delta} |\psi|^{\frac{4}{d}}\psi (s)\, ds\\
\| \psi(t)\|_{2} &\le \| \psi_{0}\|_{2}+ \int_{0}^{t} \| \psi(s)\|^{p}_{L^{2p}}\, ds,\quad p=1+\frac{4}{d}
}
By the contraction mapping principle, the Strichartz estimate 
\[
\big\| e^{-it\Delta} \psi_{0}\big \|_{L^{p}_{t}{L_{x}^{2p}}} \les \|\psi_{0}\|_{2}
\]
allows us to find a unique fixed point of the integral equation \eqref{eq:duh} in the space {   $C(\R, L^{2}(\R^{d}))\cap L^{p}_{t} L_{x}^{2p}(\R^{1+d})$} for small data.  See~\cite{Tbook, Jbook} for  introductions to the vast subject of nonlinear dispersive equations. 

\subsection{Finer restriction properties} 

The Stein-Tomas theorem~\eqref{eq:SteinTomas} is optimal for $L^{2}$ restriction. The appearance of $L^{2}$ is essential since it allows one to factor through that space in a $T^{*}T$ argument, cf.~\eqref{eq:T*T} and~\eqref{eq:dualTL}. In the applications to spectral theory this aspect is also relevant,  since it is mostly this operator which arises, rather than restriction itself. 

For the sake of completeness we nevertheless formulate the analogue of~\eqref{eq:SteinTomas} without $L^{2}$ on the left-hand side. 
The fundamental {\bf restriction conjecture} states that  in dimensions $3$ and higher 
\EQ{\label{eq:RC}
\| \widehat{f\sigma_{S}} \|_{L^{q}(\R^{d})} &\les \|f\|_{L^{\infty}(S)},\quad q>\frac{2d}{d-1}  \\
 \text{or \ \ }   &\les  \| f\|_{L^{p}(S)}, \quad p'\le \frac{(d-1)q}{d+1}
}
where {   $S\subset \R^{d}$} is the sphere or another compact surface with nonzero curvature. The range of~$q$ here is  {  \em   optimal} by 
the decay estimate 
$
|\widehat{\sigma_{S}}(\xi)|\les \langle  \xi\ra^{-\frac{d-1}{2}},
$
and the range of $p$ in the second line of~\eqref{eq:RC}. 
If true, the conjecture \eqref{eq:RC} would imply {optimal bounds} on the Hausdorff dimension of {Kakeya-Besicovitch sets} (namely that they have full dimension equal to that of the ambient space). There are other remarkable connections with number theory, and additive combinatorics. 
It is fair to say that  this conjecture, its ramifications, and other geometric/combinatorial problems connected with it such as the Erd\"{o}s distance set problem, have been the driving force behind the development of harmonic analysis over the past 20 years or so. For  many aspects of the modern theory and numerous references, see~\cite{guth}, and for a more classical survey cf.~\cite{tao}, \cite{wolff}.   
In the plane, the restriction conjecture as well as the dimension of Kakeya sets are known.

\section{Scattering theory and Fourier restriction  } \label{sec:IS}
 
We now describe a rendition of Agmon-Kato-Kuroda based on the Stein-Tomas theorem \eqref{eq:SteinTomas} rather than the trace lemma~\eqref{eq:tracelem}.  
The starting point is again the formula for the imaginary part of the resolvent, i.e., 
 the relation
 {     
 \[
 [(-\Delta - (\lambda^2 +  i0))^{-1} - (-\Delta - (\lambda^2 -  i0))^{-1}]f =  c\lambda^{-1} \widehat{\sigma_{\lambda \Sph^{d-1}}}\ast f. 
 \]
 }
Note that  the  right-hand side is precisely of the form as it appears  in the $T^{*}T$ formulation of the Stein-Tomas theorem, and thus satisfies the estimate~\eqref{eq:T*T}. 
  Kenig, Ruiz, Sogge~\cite{krs} established  the same bound for the full resolvent
   $R_{0}(\lambda)$, viz. 
\EQ{\label{eq:KRS}
\| (-\Delta -(\lambda^{2}+i0))^{-1}\|_{L^{p_{d}}(\R^{d}) \to L^{p_{d}'}(\R^{d}) } \le C \, \lambda^{-\frac{2}{d+1}}
}
As before, the question is how to transfer this result to the perturbed resolvent. To formulate the main result from~\cite{ionschlag} to this effect 
we introduce the following operators and spaces:  
\EQ{
\label{star}
M_{q}(f)(x) & : = \Big[ \int_{|y|\le 1/2} |f(x+y)|^{q}\, dy \Big]^{\frac1q}, \quad q=\max(\frac{d}{2}, 1+) \\
\| V\|_{Y} &:=  \sum_{j=0}^{\infty} 2^{j} \| V\|_{L^{\infty}(D_{j})} < \infty, \quad M_{q} V \in L^{\frac{d+1}{2}}(\R^{d}) 
}
Here $D_{j}$ are the usual dyadic shells for $j\ge1$ and $D_{0}$ is the unit ball at the origin. 
The effect of the $M_{q}$ operator is to distinguish between local singularities and decay at infinity. 
The   {    Agmon-Kato-Kuroda } theory on the basis of the Stein-Thomas type theorem~\eqref{eq:KRS} takes the following form. 

 \begin{theorem} \label{thm:IS}  Let {   $V$} be real-valued, and suppose that  {   $V=V_{1}+V_{2}$} with constituents  satisfying either of the conditions in~\eqref{star}. Then the spectrum is purely absolutely continuous, i.e.,   {   $\sigma_{ac}=[0,\infty)$}, there is no
 singular continuous spectrum, the pure point spectrum lies in {   $(-\infty,0]$}, and is discrete in {   $(-\infty,0)$}, the  eigenfunctions decay rapidly, and the wave operators {   $W_{\pm}$} exist and are complete.  
 Moreover, a suitable limiting absorption principle holds based on the spaces in~\eqref{star}. 
 \end{theorem}
 
 See the paper~\cite{ionschlag} for a precise statement of the limiting absorption principle. 
 Magnetic potentials are also admissible for this theorem, but we did not include them for the sake of simplicity.  Note that the condition  {   $M_{\frac{d}{2}} V \in L^{\frac{d+1}{2}}(\R^{d})$} is weaker in terms of decay at infinity than {   $V\in L^{\frac{d}{2}}(\R^{d})$} and {     sharp} for {   $d\ge3$}. The latter follows from an example given in~\cite{ionjer} of a potential {   $V\in L^{p}(\R^{d}), p>\frac{d+1}{2}$} with embedded eigenvalues,  and anisotropic decay
\[
 |V(x)|\simeq (1 + |x_{1}| + |x'|^{2})^{-1}. 
\]
Earlier, \cite{golsch} had established the following 
 limiting absorption principle for  {    $L^{\frac{3}{2}}$} potentials  in three dimensions:  
 
 \begin{theorem}
\label{thm:agmon} Let $V\in L^p(\R^3)\cap L^{\frac32}(\R^3), p > \frac32$
be real-valued.  Then for every $\lambda_0>0$, one has
\begin{equation}
\label{eq:V34}
\sup_{0<\eps<1,\;\lambda\ge\lambda_0}
\Big\|(-\Laplace+V - (\lambda^2+i\eps))^{-1}\Big\|_{\frac43\to 4} \le C(\lambda_0,V)\;\lambda^{-\half}.
\end{equation}
In particular, the spectrum of $-\Laplace+V$ is purely absolutely continuous on $(0,\infty)$.
\end{theorem}

Crucial to both Theorem~\ref{thm:IS} and~\ref{thm:agmon} is 
the absence of embedded eigenvalues. As discussed above, in the classical weighted $L^{2}$ context one
uses Kato's theorem for that purpose which applies to short-range potentials~\eqref{eq:short} (which is sharp in terms of point-wise decay by the famous Wigner, von Neumann potential~\cite{FrankSimon}). The results of this section, however, require a  result on the
absence of embedded eigenvalues that only assumes an
integrability condition on~$V$. One such result was obtained by
Ionescu and Jersion~\cite{ionjer}, namely:

\begin{theorem} \label{thm:IJ}
Let $V\in L^{\frac32}(\R^3)$. Suppose $u\in W^{1,2}_{\rm loc}(\R^3)$ satisfies
$(-\Laplace +V)u=\lambda^2 u$ where $\lambda\ne0$ in the sense of
distributions. If, moreover,
$\|(1+|x|)^{\delta-\half}u\|_2<\infty$ for some $\delta>0$, then $u\equiv 0$.
\end{theorem}

See also~\cite{FrankSimon}. 
The weighted $L^2$-condition with $\delta>0$ is natural in view of the Fourier
transform of the surface measure of $\Sph^2$,
which is a generalized eigenfunction
of the free case and decays like $(1+|x|)^{-1}$. Koch and Tataru~\cite{KTa} improved on this result and established absence of embedded eigenvalues assuming only $V\in L^{\frac{d+1}{2}}(\R^{d})$ as in~\eqref{star}, which is crucial for the validity of Theorem~\ref{thm:IS}. 
 
In closing let us mention other applications of Fourier restriction theory to the spectral theory of random operators. 
In \cite{bourgain} the almost sure existence and asymptotic completeness of wave operators is shown for the random lattice model $$H_{\omega}=-\Delta_{\Z^{2}}+\sum_{n\in\Z^{2}} \beta_{n}\, \omega_{n}\,\delta_{n}$$ where $\beta_{n}$ is a decaying weight $(1+|n|)^{-\frac12-\eps}$ and $\omega_{n}$ are i.i.d.~random variables such as Bernoulli. As usual, $\delta_{n}$ are the Dirac measures.  For technical reasons, Bourgain excludes energies near $0$ and near the band edges. Even though his result is formulated in the plane, the method of proof extends to all dimensions, albeit at the expense of possibly having to exclude more energies due to more complicated and singular Fermi surfaces.   The randomness therefore allows one to save half of a power of decay as compared to the short range condition~\eqref{eq:short}. 
Interestingly, this work does not invoke the Stein-Tomas theorem, but relies more on the trace lemma and entropy bounds from the probabilistic theory of Banach spaces (dual Sudakov inequality). 

In the follow up paper~\cite{bourgain2},   the Stein-Tomas theorem is used explicitly to replace the point-wise decay by an $\ell^{p}$ decay condition on the potential. He also makes an interesting reference to invoking the sharp two-dimensional restriction theory (i.e., the restriction conjecture in the plane, known as Carleson-Sj\"{o}lin or Zygmund theorem) in order to carry out a rigorous renormalization procedure. To the best of the author's knowledge this has not been carried out yet.   It remains to be seen if the modern and much more advanced  Fourier restriction theory alluded to in the previous section can be applied to spectral theory, especially in the context of the Anderson model.

\section{  Yajima's {     $L^{p}$} theory for the intertwining operator    }

In the 1990s {    Kenji Yajima} initiated a far-reaching investigation of the $L^{p}$ boundedness properties of the wave operators~\eqref{eq:W}. His starting point was the stationary representation of the wave operators due to Kato~\cite{kato}.  We cannot give a complete account here of the results~\cite{yajima0}--\cite{yajima8}, \cite{AY}, \cite{Wed} but instead state some representative theorem: one has the boundedness {   $W_{\pm} : L^{p}(\R^{3})\to L^{p}(\R^{3})$, for all $1\le p\le \infty$,  provided    $|V(x)|\le \langle x\ra^{-5-\eps}$  and provided there is   no zero energy eigenvalue or resonance.  The latter assumption is essential as we will see in the next paragraph. In fact, if    zero energy is singular, then the wave operators are $L^{p}(\R^{3})$ bounded only in the smaller range   $3/2<p<3$, and provided one has the decay $|V(x)|\le \langle x\ra^{-6-\eps}$.  

Similar results hold in higher dimensions, but not only is more decay required, but some regularity on the potential is needed as well. 
Low dimensions behave differently, and boundedness at $p=1,\infty$ is lost for the line and the plane. This has to do with the appearance of a Hilbert transform in the kernel representation of the wave operators.

\smallskip  

Apart from its intrinsic interest in terms of shedding more light on the nature of the wave operators, Yajima's theory provides a very quick way of obtaining dispersive estimates for operators with  a potential from the free case. To be more specific, consider the evolution 
  $e^{it \Phi(H)} P_{c}(H)$} where $\Phi$ is a polynomial, say or some other function taking $H$ to the self-adjoint operator $\Phi(H)$, and $P_{c}(H)$ is the projection onto the absolutely continuous spectrum (we have asymptotic completeness of $H$). As usual, $H_{0}=-\Delta$. Then  
\[
e^{it \Phi(H)} P_{c}(H) = W e^{it \Phi(H_{0})} W^{*}
\]
allows one to transfer $L^{p}$ or Strichartz estimates from $H_{0}$ to $H$ provided $0$ energy is regular, simply by bounding $W$ and $W^{*}$ by their $L^{p}$ operator norms. The 
importance of the {   $0$} energy condition is implied by this, too.  For example,  in three dimensions one has 
\EQ{\label{eq:rauch}
\big\| e^{it H} f \big\|_{\infty } \le \|W\|_{\infty\to\infty}\|W\|_{1\to1}\; C t^{-\frac32}\|f\|_{1},\qquad f\perp \text{bound states}
}
This is known to fail in the presence of a $0$ energy obstruction. In fact, a whole power of $t$ is lost from the decay in that case, cf.~\cite{jenkat}, \cite{rauch}, \cite{mur1,mur2}, \cite{JSS}.  

In some applications it might not be possibly to invoke the $L^{p}$ theory of $W_{\pm}$ to two reasons: 
(i) the assumptions on potential are too strong (ii) in some nonlinear applications $0$ energy singularities do arise. 

Both of these issues occur for example in Krieger's work with the author~\cite{KrSc}. This work deals with the conditional (in the spirit of a center-stable manifold) asymptotic stability analysis of an unstable soliton for the energy critical radial nonlinear wave equation
\[
u_{tt} -\Delta u - u^{5} =0
\]
in $\R_{t,x}^{1+3}$. The unique radial stationary $\dot H^{1}(\R^{3})$ solution to this equation is the Aubin-Talenti solution $W(x) = (1+|x|^{2}/3)^{-\frac12}$ and its rescalings $W_{\lambda}(x) = \lambda^{\frac12}W(\lambda x)$.    The linearization about $W$ leads to the operator 
\[
H = -\Delta - 5W^{4}(x)
\]
which has the $0$ energy resonance given by $\psi(x)= \partial_{\lambda}\Big|_{\lambda=1}W_{\lambda}(x)$. Note that $H\psi=0$ and $\psi(x)\sim |x|^{-1}$ as $x\to\infty$. So $\psi$ is indeed not a $0$ energy eigenfunction, but rather  a resonance.  In addition, $H$ has a unique negative eigenvalue. 

It is shown in~\cite{KrSc} that the wave evolution $\frac{\sin(\sqrt{H} t)}{\sqrt{H}}P_{c}(H)$ does not decay at the free rate $t^{-1}$; rather, along the resonance it does not decay at all. However, after subtracting the rank~$1$ operator $\psi\otimes\psi$, the free rate of $t^{-1}$ is regained. Such results currently fall outside of the scope of the wave operator results, and require a direct approach based on the Laurent expansion of the resolvent near~$0$. On the other hand, energies bounded away from~$0$ are not so much the issue here, and can be dealt with in a more general fashion. 

\subsection{  Yajima's proof, expansion of the wave operators    }

In the remainder of this section we present some of the essential steps in Yajima's analysis of $W_{\pm}$ in three dimensions. 
Relying on the time-dependent representation of the wave operators as in~\eqref{eq:cook}, we 
iterate the Duhamel formula to  obtain the expansion,  formally at first,
\EQ{\label{eq:dyson}
W f &= f + W_{1} f + \ldots + W_{n} f + \ldots,\\
\nonumber W_{1} f &= i \int_{t>0} e^{-i t \Delta} V e^{i t \Delta} f \dd t,\ \ldots \\
W_{n} f &= i^n \int_{t>s_1>\ldots>s_{n-1}>0} e^{-i(t-s_1)\Delta} V e^{-i(s_1-s_2) \Delta} V \ldots \\
\nonumber &\qquad\qquad e^{-i s_{n-1} \Delta} V e^{it\Delta} f \dd t \dd s_1 \ldots \dd s_{n-1}
}
  for all  $f\in L^2$.    There are a number of ways in which one might justify this series expansion rigorously (which is nothing other than the Dyson series). One elegant, and essentially optimal way in terms of the conditions on~$V$, would be to use the 
 Keel-Tao Strichartz endpoint (in {   $\R^{3}$}) which is of the form
{   \EQ{\nn 
\|e^{it H_0} f\|_{L^2_t L^{6, 2}_x} &\les \|f\|_{L^2}\\
\Big\|\int_\R e^{-isH_0} F(s) \dd s\Big\|_{L^2_x} &\les\|F\|_{L^2_t L^{6/5, 2}_x}, 
}
see~\cite{keetao}.  The spaces $L^{p,q}$ are the Lorentz spaces, see for example~\cite{bergh}. The operator acting by multiplication by the potential $V$ satisfies 
$$V:L^{6, 2}_x(\R^{3})\to L^{6/5, 2}_x(\R^{3})$$  provided $V\in L^{\frac32,\infty}(\R^{3})$ (the weak $L^{\frac32}(\R^{3})$ space). 
}
We thus conclude that the Dyson series converges in {   $L^{2}(\R^{3})$}  if the potential is small {   $\| V\|_{3/2,\infty}\ll 1$}.  We remark that the $L^{\frac32}(\R^{3})$  condition on~$V$ enjoys a special significance. It is precisely the scaling invariant space associated with the Schr\"{o}dinger operator, which has the underlying scaling $V\to \lambda^{2}V(\lambda x)$. This is valid in all dimensions, and the invariant norm under this symmetry is~$L^{\frac{d}{2}}(\R^{d})$.  In terms of power law decay, the critical (scaling invariant) fall off rate is~$|x|^{-2}$. The importance of this threshold is well-known, see for example~\cite{BPST1, BPST2} in the context of linear dispersive estimates.  It is important to note the difference between the $|x|^{-1}$ threshold in the short range condition~\eqref{eq:short} and the aforementioned $|x|^{-2}$ decay which is critical relative to scaling. Loosely speaking, the former is the natural decay rate for the {\em scattering theory} which hinges on  the resolvent for  large or at least positive energies, whereas the latter is the cutoff for the dispersive theory for which the $0$ energy behavior of the resolvent is the deciding factor. 

If $V$ is not small, then the infinite series is much more delicate. One option is to truncate, i.e., 
\EQ{\label{eq:trunc}
W f = f + W_{1} f + \ldots + W_{n-1} f + \tilde W_{n} f
}
where 
\EQ{\nn 
\tilde W_{n} f &= i^n \int_{t>s_1>\ldots>s_{n-1}>0} e^{i(t-s_1)(-\Delta+V)} V e^{-i(s_1-s_2) \Delta} V \ldots \\
\nonumber &\qquad\qquad e^{-i s_{n-1} \Delta} V e^{it\Delta} f \dd t \dd s_1 \ldots \dd s_{n-1}
}
In other words, the perturbed evolution remains in the Duhamel formula (but it need not appear in the first position, and can be moved to the right).  This is the approach chosen in~\cite{yajima0}. It leads to losses in terms of the decay of~$V$.  In the following section we will present Wiener's theorem in convolution algebras as a summation method for this infinite series. 

\subsection{ Representations of the summands $W_{n}$     }

Following Yajima, we now exhibit the structure of the individual terms $W_{n}$ in~\eqref{eq:dyson} in order to show how $L^{p}$ boundedness arises.  
First we introduce convergence factors $e^{-\eps t}$ in the Duhamel expansion, i.e., we  introduce the regularized operators
\EQ{ \label{eq2.5}
W_{n+}^\eps  f&:= i^n \int_{0\leq t_1 \leq \ldots \leq t_n} e^{i(t_n-t_{n-1})H_0-\eps (t_n-t_{n-1})} V \ldots \\
& e^{i(t_2 - t_1)H_0 - \eps (t_2 - t_1)} V e^{it_1 H_0 - \eps  t_1} V e^{-it_nH_0} f \dd t_1 \ldots \dd t_n,
}
together with
\EQ{\label{eq2.6}
W_+^\eps  = I + i \int_0^{\infty} e^{it H-\eps  t} V e^{-it H_0} \dd t.
}
for $\eps>0$. 
Taking Fourier transforms then yields, for $V, f,g$ Schwartz functions: 
\EQ{\label{eq:keyrep1}
& \langle W_{n}^\eps  f, g\rangle = \\
&  \frac {(-1)^n}{(2\pi)^3} \int_{\R^{3(n+1)}} \!\!\!\!\! \frac{\prod_{\ell=1}^n \widehat V(\xi_{\ell} - \xi_{\ell-1}) \dd \xi_1 \ldots \dd \xi_{n-1}}{\prod_{\ell=1}^n (|\eta+\xi_{\ell}|^2- |\eta|^2 + i \eps )} \widehat f(\eta) \overline{\widehat g}(\eta+\xi_n) \dd \eta \dd \xi_n 
}
as well as 
\EQ{\label{eq:keyrep2}
& \langle W_{1+}^\eps  f, g\rangle = - \frac 1 {(2\pi)^3} \int_{\R^6} \frac {\widehat V(\xi)}{|\eta+\xi|^2- |\eta|^2 + i \eps } \widehat f(\eta) \ov{\widehat g}(\eta+\xi) \dd \eta \dd \xi \\
& \qquad\qquad = \int_{\R^{6}} K_{1}^{\eps}(x,x-y) f(y)\, dy \, \ov g(x)\, dx 
}
The reason for writing $K(x,x-y)$ rather than $K(x,y)$ in the previous line lies with the fact that we obtain a cleaner expression for this kernel. In fact, one has 
\EQ{\nn 
K^\eps _{1}(x,z) &= c |z|^{-2}  \int_0^\infty e^{-is\hat{z}\cdot(x-z/2)} \widehat{V}(-s\hat{z}) e^{-\eps \frac{|z|}{2s}}\: s\,ds,\quad \hat{z}=z/|z|\\
 K_{1}(x,z) &= c|z|^{-2}L(|z|-2x\cdot\hat{z},\hat{z}), \;\; L(r,\omega)= \int_{0}^{\infty} \widehat{V}(-s\hat{z}) e^{i\frac{rs}{2}} \,s\, ds
}
where we passed to the limit $\eps\to0$ in the last line.   The details of these computations can be found in~\cite{becsch1}. 

\subsection{  The structure of  $W_{1}$ in $\R^{3}$    }\label{sec:W1struct}

Denote by 
 {   $S_\omega x := x-2(\omega\cdot x) \omega$} the  reflection about the plane   $\omega^\perp$.  In view of the preceding, 
\EQ{\nn 
 (W_{1}f)(x) & = \int_0^\infty \int_{\Sph^2 } L(r-2\omega\cdot x,\omega) f(x-r\omega)\, drd\omega \\
&= \int_{\Sph^2 } \int_{\R} \one_{[r>-2\omega\cdot x]} L(r,\omega) f(S_{\omega }x -r\omega)\, drd\omega \\
& = \int_{\Sph^2 }\int_{\R^3} g_1(x,dy,\omega) f(S_\omega x-y)\,  d\omega
}
Therefore,  with {   $\calH^1_{\ell_\omega}$} the Hausdorff measure on the line along    $\omega$,
\EQ{\label{eq:L}
 g_{1}(x,dy,\omega) &:= \one_{[(y+2x)\cdot\omega>0]} L(y\cdot\omega,\omega) \, \calH^1_{\ell_\omega}(dy)  \\
  \int_{\Sph^2 }  \| g_1(x,dy,\omega)\|_{\mes_{y} L^\infty_x}  \, d\omega &\le \int_{\Sph^2 }\int_{\R} |L(r,\omega)|\, drd\omega  =: \|L\| \\
 \| W_{1}f\|_{p}  &\le \|L\| \|f\|_{p}
}
This shows that (i) $W_{1}$ acts as a weighted average of translations and reflections (ii) $W_{1}$ is therefore $L^{p}$ bounded uniformly in $p$ provided $L$ has  finite norm as above.

\subsection{  Bounding $L$     }

We now address the boundedness of the function $L$, cf.~\eqref{eq:L}. 
Define
{   \[
\|f\|_{B^{\beta}}:= \| \one_{[|x|\le 1]} f\|_{2} + \sum_{j=0}^{\infty} 2^{j\beta} \big\| \one_{[2^{j}\le |x|\le 2^{j+1}]} f \big\|_{2} <\infty
\] }
Then {   $\dot B^{\frac12} \embed  L^{\frac32,1}(\R^3)$, $\dot B^{1} \embed  L^{\frac65,1}(\R^3)$}, and 
\EQ{ \nn 
\|L(r,\omega)\|_{L^2_{r, \omega}} & \les \|V\|_{L^2} \\
\|L(r,\omega)\|_{L^1_{r, \omega}}  & \les \sum_{k \in \Z} 2^{k/2} \|\one_{[2^k, 2^{k+1}]}(|r|) L(r,\omega)\|_{L^2_{r, \omega}} \les \|V\|_{\dot B^{\frac12}}\les \|V\|_{B^{\frac12}}
} 
Hence, $L^{p}$ boundedness of $W_{1}$ holds under the scaling invariant condition $\|V\|_{\dot B^{\frac12}}<\infty$ (in particular, it is enough if $\|V\|_{L^{\frac32,1}(\R^{3})}<\infty$). As far as the higher terms $W_{2}$ etc.\ are concerned, 
Yajima showed that for small potentials in the sense {   $\| V\|_{B^{1+\eps}}\ll 1$} one has 
{   
\[
\| W_{n} f\|_{p}\le C^{n} \| V\|_{B^{1+\eps}}^{n}\, \|f\|_{p}
\]
}
Note that this is off the scaling critical value by $\frac12+\eps$. 
 The Dyson series \eqref{eq:dyson} can thus be summed in the $L^{p}$ operator norm leading to the full result for small potentials in $B^{1+\eps}(\R^{3})$. 
For large potentials losses appear by terminating the expansion as in~\eqref{eq:trunc}, which is why one ends up with the stronger decay requirement $|V(x)|\les \langle x\rangle^{-6-\eps}$.  This is due to the presence of the perturbed evolution $e^{itH}$
in the final term $\tilde{W_{n}}$. For large enough~$n$, the many free evolutions appearing in that term have a regularizing effect. 

\section{ Structure Theorems    }\label{sec:wiener}

In this section we discuss the following two results from \cite{becsch1,becsch2}.  They establish that the full wave operator retains a structure similar to that of~$W_{1}$ above for large potentials, albeit under a stronger condition on~$V$ than we needed for $W_{1}$. 

\begin{theorem}[Beceanu-S. 16] \label{thm:main1}
Let $V \in B^{1+}$  real-valued, and assume that zero energy is regular for $H = -\Delta+V$ (i.e., no eigenvalue or resonance at~$0$). There exists  $g(x, dy, \omega) \in L^1_\omega \mc M_y L^\infty_x$ with 
\EQ{\nn 
 & \int_{\Sph^2} \|g(x, dy, \omega)\|_{\mc M_y L^{\infty}_x}  \dd \omega < \infty \\
(W_{+} f)(x) &= f(x) + \int_{\Sph^2} \int_{\R^3} g(x, dy, \omega) f(S_\omega x - y)   \dd \omega.
} 
Suppose $X$ is a Banach space of measurable functions on $\R^3$, which is invariant under translations and reflections, and so that Schwartz functions are dense (or dense in $Y$ with  $X=Y^*$). 
Assume  $\|\one_H  f\|_X\le A\|f\|_X$ for all half spaces $H\subset\R^3$ and $f\in X$ with some uniform constant~$A$. Then 
\EQ{\nn 
\| W_+ f\|_X \le AC(V) \|f\|_X \qquad\forall\; f\in X
}
where $C(V)$ is a constant depending on $V$ alone. 
\end{theorem}

In particular, one has $L^{p}(\R^{3})$ boundedness uniformly in $1\le p\le \infty$. This theorem is sharp in several ways. On the one hand, it cannot hold in dimension $d=1$ since $L^{p}$ boundedness fails at $p=1$ and $p=\infty$ by~\cite{Wed}. It presumable also fails in the plane. On the other hand,  the assumption that~$0$ energy is regular is essential, too. Otherwise we might deduce dispersive decay for $e^{itH}P_{c}$ which is known to fail, see the discussion around~\eqref{eq:rauch}.  A quantitative version of the previous theorem is also available.

\begin{theorem}[Beceanu-S. 17]  \label{thm:main2}
Let $V\in B^{1+2\gamma}$, $0<\gamma$, with the same $0$ energy  hypothesis as above.
Then 
\EQ{\label{eq:M0}
& \int_{\Sph^2 }  \|  g(x,dy,\omega)\|_{\mes_{y} L^\infty_x}  \, d\omega  \le C_{0} (1+\|V\|_{B^{1+2\gamma}})^{38+\frac{105}{\gamma}} (1+M_{0})^{4+\frac{3}{\gamma}} \\
& \sup_{\eta\in\R^{3}} \sup_{\eps>0}\big \|   \big(I + R_0(|\eta|^2 \pm i\eps ) V  \big)^{-1} \big\|_{\infty \to \infty} =: M_{0}<\infty
} 
where  $C_{0}$  absolute constant. 
\end{theorem} 

We make the following remarks: 

\begin{itemize}
\item {   $0$} energy regular means exactly that {   $$ M_{00}:=\|   \big(I + (-\Delta)^{-1} V  \big)^{-1} \big\|_{\infty \to \infty}  <\infty.$$}  It is shown in~\cite{becsch1} that this is equivalent to the definition by~\cite{jenkat}. Moreover, it is established in loc.~cit.\  that under this assumption $M_{0}<\infty$. This relies on the results of Section~\ref{sec:IS}, in particular on Theorem~\ref{thm:IS} and thus depends on  Stein-Tomas type Fourier restriction bounds. 
\item It would be desirable to bound {   $M_{0}$} through $M_{00}$ and the size of $V$ in some sense. The control of {   $M_{0}$} is not effective, see the discussion above about the limiting absorption principle, and the effective estimates of~\cite{RT}. 
\item These results fall short by more than {    $\frac12$} from the scaling invariant class {    $\dot B^{\frac12}$}. It is not clear how to avoid this loss within the framework of~\cite{becsch1}. It is conceivable that the methods are optimal {\em assuming only decay of} $V$, and thus a scaling invariant condition on $V$ leading to a structure result as above would need to involve spaces which measure more than decay in $L^{2}$ along dyadic shells. 
At the end of this note we state the result from~\cite{becsch2} which builds such a space, but only for small $V$. 
\item It is likely that Theorem~\ref{thm:main1} remains valid in $B^{1}$, but no quantitative analogue of it as in Theorem~\ref{thm:main2} would then hold (within the confines of the methods of \cite{becsch1, becsch2}). 
\item A version of Theorem~\ref{thm:main1} should hold in higher dimensions. 
 \end{itemize}

For the remainder of this paper we give some indications of the methods involved in proving these theorems, in particular, the Wiener inversion technique. 

\subsection{  Wiener algebra and inversion    }

We cannot sum the Dyson series~\eqref{eq:dyson}. Instead we use {    Beceanu's operator-valued Wiener formalism}, cf.~\cite{bec,becgol}. First recall the classical Wiener theorem for the convolution algebra on the line.  

\begin{prop}
Let $f\in L^1(\R^{d})$. There exists $g\in L^1(\R^{d})$ with
\EQ{
(1+\hat{f}\,)(1+\hat{g})=1 \text{\ \ on\ \ }\R^{d}
}
iff $1+\hat{f}\ne0$ everywhere. 
Equivalently, there exists $g\in L^1(\R^{d})$ so that 
\EQ{
(\delta_0+f)\ast (\delta_0+g)=\delta_0
}
 iff $1+\hat{f}\ne0$ everywhere on $\R^{d}$.  The function $g$ is  unique. 
\end{prop}

There are two critical features in the classical proof (which is presented in~\cite{becsch1}): 
\begin{itemize}
\item each $f\in L^{1}(\R)$ exhibits a uniform {   $L^{1}$}-modulus of continuity under translation.
\item one has vanishing at {   $\infty$} in the  {   $L^{1}$} sense (``no mass at infinity''). 
\end{itemize}
There is a well-known alternative proof via Gelfand-Naimark theory, which depends on identifying the maximal ideal space. We do not follow that approach here since it appears to work only in the Abelian setting. In our spectral theory setting, however, the algebras are noncommutative. We now present such an algebra.

\subsection{  An operator-valued version    }\label{sec:52}

Let {   $X$} be a Banach space,  and denote by {   $\W_{X}$}  the  algebra of bounded linear
maps {   $T: X \to L^1(\R; X)$} with the convolution 
{   \[
S \ast T(\rho)f = \int_\R S(\rho-\sigma)T(\sigma)f\, d\sigma
\]
}
As usual, one adjoins unit $\delta$, and we  denote this  larger algebra by {   $\widetilde \W_{X}$}. 
The Fourier transform exists and  satisfies
{   
$$\sup_{\lambda} \norm[\hat{T}(\lambda)][\B(X)] \le \norm[T][\W_X]$$
}
In this setting one has the following exact analogue of the scalar Wiener theorem. Note how the two conditions below 
capture the continuity relative to translations, and the vanishing at~$\infty$. 

\begin{theorem}[\cite{bec},\cite{becgol}] \label{thm:BG}
 Suppose $T\in \W_X$ satisfies 
\begin{enumerate}
\item  $\lim\limits_{\delta \to 0} 
\norm[T(\rho) - T(\rho-\delta)][\W_X] = 0$.
\item  $\lim\limits_{R \to \infty}
\norm[T\chi_{|\rho| \ge R}][\W_X] = 0$.
\end{enumerate}
If $I + \hat{T}(\lambda)$  is invertible in $\B(X)$ for all
$\lambda$, then ${\mathbf 1} + T$ possesses an inverse in 
$\widetilde{\W}_X$ of the form ${\mathbf 1} + S$.
\end{theorem}

\subsection{ Application to dispersive estimates     } \label{sec:BGdisp}

We now apply Theorem~\ref{thm:BG} to derive decay of the  Schr\"{o}dinger evolution in~$\R^{3}$.

Set    $$R_{0}^{-}(\lambda^{2})(x)=(4\pi|x|)^{-1}e^{-i\lambda|x|}, \; \widehat{T^{-}}(\lambda)=  V R_{0}^{-}(\lambda^{2}).$$  Then 
{   \begin{equation} \label{eq:Tminus}
T^-(\rho)f(x) = (4\pi \rho)^{-1} V(x) \int_{|x-y|= \rho} f(y)\,dy
\end{equation}}
and thus 
{   \begin{align*}
\int_{\R^3}\int_\R |T^-(\rho)f(x)|\, dx\,d\rho &\le
\frac{1}{4\pi} \int_{\R^3} \int_{\R^3} \frac{|V(x)|}{|x-y|}|f(y)|\,dy\,dx \\
&\le \frac{1}{4\pi} \norm[V][\Kato] \norm[f][1].
\end{align*}
}
where {   $\norm[V][\Kato] = \| |x|^{-1}\ast |V| \|_{\infty}$.} Hence our underlying algebra is {   $\W_{L^{1}}$}. The crucial  pointwise invertibility condition  needed in Wiener's theorem takes the following form:
 {   $$(I+ V R_{0}^{-}(\lambda^{2}))^{-1}\in \B(L^{1})$$}
 which is precisely the invertibility of the Birman-Schwinger operator relative to $L^{1}$. 
By spectral and scattering theory this is guaranteed for positive energies, and for zero energy it becomes an assumption. 
The invertibility of $I+ T^-(\rho)$ therefore holds in~$\widetilde{\W}_{L^{1}}$. 
 This is used in~\cite{becgol}  to  prove dispersive estimates for Schr\"odinger in {   $\R^{3}$} for {   $\|V\|_{\K}<\infty$}. 

To place this into context, my earlier~\cite{rodsch}  had shown that if 
\EQ{\label{eq:klein}
\sup_{x\in\R^{3}} \int_{\R^{3}} \frac{|V(y)|}{|x-y|}\, dy <4\pi
}
then for {   $V$} real valued one has dispersive estimate 
\EQ{\label{eq:IRklein}
\| e^{itH} f\|_{\infty} \le C|t|^{-\frac32} \|f\|_{1},\quad H=-\Delta + V
}
The strategy was to write the evolution via the functional calculus, invoke that the density of the spectral measure is the imaginary part of the resolvent,  and to expand the resolvent into an infinite Born  series. Each term in this series is then handled by certain  oscillatory integral estimates.  The series converges because of the smallness condition~\eqref{eq:klein}. 

It remained an 
{    open problem} to obtain the analogue of~\eqref{eq:IRklein} (on the orthogonal complement of the bound states) for large~$V$. 
This problem was solved in~\cite{becgol} assuming zero energy is regular (which is a necessary condition) by means of this  Wiener algebra.

\subsection{  Algebra for intertwining operators    }

The formalism underlying Theorems~\ref{thm:main1} and~\ref{thm:main2} is much heavier than the one in Section~\ref{sec:52}. 
The formulas for {   $W_{n}$} suggest using three-variable kernels. To begin with, we introduce an algebra which will be far too weak to control the Dyson series~\eqref{eq:dyson}. But it is indispensable as an ambient space which will contain the key algebra $Y$, see below. Define the following space as a subset of tempered distributions 
{   
\EQ{ \nn 
Z & := \{T(x_0, x_1, y)\in \calS'(\R^{9}) \mid \;   \calF_{y} T(x_{0},x_{1},\eta) \in L^{\infty}_{\eta}L^{\infty}_{x_{1}}L^{1}_{x_{0}} \}  \\
\| T\|_Z & := \sup_{\eta\in\R^3} \| \mc F_{y}T(x_0, x_1, \eta)\|_{{L}^{\infty}_{x_1}L^{1}_{x_{0}}} 
} }
The convolution operation {    $\oast$}  on    $T_{1},T_{2}\in Z$  is 
\EQ{\nn 
(T_{1}\oast T_{2} )(x_0, x_2,y) = \mc F_{\eta}^{-1} \Big[ \int_{\R^{3}}\mc F_{y}T_{1}(x_0, x_1, \eta) \mc F_{y}T_{2}(x_1, x_2, \eta)\, dx_{1} \Big](y)
}
where the Fourier transform is understood in the sense of tempered distributions. 
Furthermore, introduce the seminormed space   $V^{-1}B$  defined as 
$$
V^{-1}B=\{f\ \text{measurable} \mid V(x) f(x) \in B^{\sigma} \}
$$ 
with the seminorm {    $\|f\|_{V^{-1} B}:=\|V f\|_{B^{\sigma} }$}.   
Finally,  set {    $X_{x,y}:=L^{1}_{y}V^{-1}B_{x}$}.  Then {   $L^{1}_{y}L^{\infty}_{x}$} dense in {   $X_{x,y}$}. 

Now introduce the following key space {    $Y$}   of three-variable kernels   
\EQ{\nn 
Y &:= \Big\{T(x_0, x_1, y) \in Z \mid \; \forall f \in  L^\infty\\
& (fT)(x_{1}, y):= \int_{\R^3} f(x_0) T(x_0, x_1, y) \dd x_0 \in X_{x_{1},y} \Big\},
}
with norm 
\EQ{\nn 
\|T\|_Y &:= \|T\|_Z  + \|T\|_{\B(V^{-1}B_{x_0}, X_{x_{1},y} )} 
}
For {   $\frak X \in L^{1}_{y}L^{\infty}_{x}$}, define the operation of  \emph{contraction} of {   $T\in Y$} by {   $\frak X$} to be  
 \[
(\frak X T)(x, y) := \int_{\R^6} \frak X(x_0, y_0) T(x_0, x, y-y_0) \dd x_0 \dd y_0.   
\] 
Then   $\frak X T\in X_{x,y}$,  and $\|\frak X T\|_{X}\le \|T\|_{Y}\|\frak X\|_{X}$. This property turns {   $Y$} into an algebra under {   $\oast$}.  

To understand the reason behind these structures, we return to the explicit formulas~\eqref{eq:keyrep1}, \eqref{eq:keyrep2} obtained above for $W_{n}$. 
In order to relate them to the operation of convolution~$\oast$, we 
 define 
  \EQ{ \label{eq:T1}
\mc F_y T_{1+}^\eps (x_0, x_1, \eta)  
&=  e^{-ix_1 \eta}\,  R_0(|\eta|^2-i\eps )(x_0, x_1)V(x_0)\,  e^{ix_0 \eta} \\
T_{2+}^\eps &= T_{1+}^\eps\oast T_{1+}^\eps,\; T_{3+}^\eps = T_{2+}^\eps\oast T_{1+}^\eps  \text{\ \ etc.} 
} 
Then  
\EQ{\label{eq:Wneps}
\langle W_{n+}^\eps  f, g \rangle &= \frac{(-1)^n}{(2\pi)^3} \int_{\R^6} \mc F^{-1}_{x_0} \mc F_{x_n, y} T_{n+}^\eps (0, \xi_n, \eta) \widehat f(\eta) \ov {\widehat g}(\eta + \xi_n) \dd \eta \dd \xi_n \\
&= (-1)^n \int_{\R^9} \mc F_{x_{0}}^{-1 }T_{n+}^\eps (0, x, y) f(x-y) \ov g(x) \dd y \dd x.
}
Replacing the free resolvent in~\eqref{eq:T1} with the perturbed one yields 
 $T_\pm^\eps $ which is  given by the distributional Fourier transform 
\EQ{\label{eqn2.13}
\mc F_{y} T_\pm^\eps (x_0, x_1, \eta) &:= e^{ix_0 \eta} \big(R_V(|\eta|^2 \mp i\eps ) V\big)(x_0, x_1) e^{-ix_1\eta};
}
where we assume that $0$ energy is regular for $H=-\Delta+V$. 
In view of~\eqref{eq2.6} we conclude in analogy with~\eqref{eq:Wneps} that 
\EQ{\nn 
& \langle W_{+}^\eps  f, g \rangle \\
&= \langle f,g\ra - \frac{1}{(2\pi)^3} \int_{\R^6} \mc F^{-1}_{x_0} \mc F_{x_1, y} T_{+}^\eps (0, \xi_1, \eta) \widehat f(\eta) \ov {\widehat g}(\eta + \xi_1) \dd \eta \dd \xi_n \\
&= \langle f,g\ra  -  \int_{\R^9} \mc F_{x_{0}}^{-1 }T_{+}^\eps (0, x, y) f(x-y) \ov g(x) \dd y \dd x.
}
Note the similarity of the right-hand sides of \eqref{eq:T1} and \eqref{eqn2.13} with the Birman-Schwinger operator $R_{0} V$ which plays a  prominent role  in scattering theory and the Beceanu-Goldberg result of Section~\ref{sec:BGdisp}. The difference here is that the operators appearing in \eqref{eq:T1} and \eqref{eqn2.13}  are of this type, but conjugated by the modulation operator $(M_{\eta}f)(x) = e^{i\eta x} f(x)$. In contrast to the algebras in Section~\ref{sec:BGdisp} the energy parameter in this context is {\em truly three dimensional} precisely because of these modulations, whereas inside the resolvent it only appears through its length~$|\eta|$. This is also the reason why our algebras are considerably more complicated as compared to the previous section.

\subsection{  Key invertibility problem     }

In classical scattering theory, the perturbed resolvent $R_{V}(z)$ is controlled from the resolvent identity~\eqref{eq:RI} and the inversion of the Birman-Schwinger operator $I+R_{0}(z)V$ (in~\eqref{eq:symmRI} we chose a symmetric form for convenience only,  since we then can invert in~$L^{2}$ rather than in weighted $L^{2}$).  There is a completely analogous inversion problem in the algebra~$Y$ that we face here.  
To begin with, due to the fact that the aforementioned modulations $M_{\eta}$ cancel each other under operator composition, we note that analogue of the
resolvent identity for the operators  
   $T_{1+}^{\eps}$ and $T_{+}^{\eps}\in Z$  reads as follows:  
\EQ{\label{keyinv}
(I+T_{1+}^{\eps})\oast (I-T_{+}^{\eps})=  (I-T_{+}^{\eps}) \oast (I+T_{1+}^{\eps}) = I
}
This identity is valid in the ambient algebra~$Z$. The key invertibility problem that we face in establishing the structure formulas for $W_{\pm}$ is 
to show that we may solve~\eqref{keyinv} for $T_{+}^{\eps}$ in the much smaller 
 algebra   $Y$.  To phrase this differently: 
{\em 
If {   $I + T_{1+}^\eps $ } is invertible in {   $Y$}, hence in {   $Z$}, its inverse is {   $I - T_+^\eps $} both in {   $Z$} and in {   $Y$}, 
hence we obtain that {   $T_+^\eps  \in Y$} uniformly in {   $\eps>0$}.
}

\smallskip 

As a first step, one needs to address that $T_{1+}^{\eps}$ belongs to $Y$. This is done in Lemma~6.2 and Corollary~7.4 in~\cite{becsch1}. To summarize what is done there, define {   $Y$} with {   $\sigma\ge \frac12$} fixed. Then 
\EQ{\label{eq:T1Y} 
\sup_{\eps>0} \| T_{1+}^{\eps} \|_{Y} &\les \| V\|_{B^{\frac12+\sigma}} \text{\ \ whence by induction}\\
\sup_{\eps>0} \| T_{n+}^{\eps} \|_{Y} &\le C^{n} \| V\|_{B^{\frac12+\sigma}}^{n} \text{\ \ for all \ \ } n\ge 1
}
Note that due to $\sigma\ge\frac12$ one loses $\frac12$ of a power of decay which is reflected in Theorem~\ref{thm:main1}. It is not clear how to avoid this loss in this exact framework, and possibly~\eqref{eq:T1Y} is optimal. 

\subsection{  Recursive definition of the structure functions  for $W_{n}$  }

To illustrate the usefulness of this algebra formalism, we now very easily obtain a structure formula for $W_{n}$ in analogy to the one derived for $W_{1}$ in Section~\ref{sec:W1struct}. In fact, we claim that 
\EQ{\label{eq:Wng}
(W_{n+}f)(x) = \int_{\Sph^2 }\int_{\R^3} g_n^{\eps}(x,dy,\omega) f(S_\omega x-y)\,   d\omega 
}
where  for fixed {   $x\in\R^3$, $\omega\in \Sph^2 $}  the expression {   $g_n^{\eps}(x,\cdot,\omega)$} is a measure satisfying 
 \EQ{\nn
\sup_{\eps>0} \int_{\Sph^2 } \| g_n^{\eps}(x,dy,\omega)\|_{\mes_{y} L^\infty_x}  \, d\omega \le  C^{n}\|V\|_{B^{\frac12+\sigma}}^{n}
}
Identifying the operator {   $W_{n+}^{\eps}$} with its kernel one has 
  \EQ{\nn
W_{n+}^{\eps} &= (-1)^{n} \one_{\R^{3}}T_{n+}^{\eps} = (-1)^{n} \one_{\R^{3}}(T_{(n-1)+}^{\eps}\oast T_{1+}^{\eps} )\\
&= -((-1)^{n-1} \one_{\R^{3}} T_{(n-1)+}^{\eps}) T_{1+}^{\eps} = - W_{(n-1)+}^{\eps} T_{1+}^{\eps}
}
In the second line we are  contracting  a kernel in    $Y$ by an element of  $X$. Thus
\EQ{\label{eq:WnX}
\sup_{\eps>0} \| W_{n+}^{\eps}\|_{X}\le \|\one_{\R^{3}}\|_{V^{-1}B} \sup_{\eps>0} \|T_{n+}^{\eps}\|_{Y} \le C^{n}\|V\|_{B^{\frac12+\sigma}}^{n+1}
}
and with {   $f^{\eps}_{y'}(x')=W^{\eps}_{(n-1)+}(x',y')$ } we have  
\EQ{\label{eq:key gn}
g_{n}^{\eps}(x,dy, \omega) := \int_{\R^{3}} g_{1,f^{\eps}_{y'}}^{\eps} (x,d(y-y'),\omega) \, dy' \\
}
where {   $g_{1,f^{\eps}_{y'}}^{\eps}$} is the structure function for $W_{1}$ associated with the potential~$f^{\eps}_{y'} V$.  See Proposition~7.6 in~\cite{becsch1} for more details. 

\subsection{  Wiener theorem in $Y$    }

For {\em small} potentials we can sum the structure formulas~\eqref{eq:Wng} and obtain the structure formula for $W_{\pm}$ as in Theorem~\ref{thm:main1}.  For large potentials we need to resort again to a Wiener formalism in order to solve equation~\eqref{keyinv} above. 
The precise formulation of the Wiener theorem for the algebra~$Y$ which was used in~\cite{becsch1} reads as follows. The space $\mc F Y$ refers to the Fourier transform of~$Y$ relative to the $y$ variable. 

\begin{theorem}
Suppose $V\in B^{\sigma}$ with $\frac12\le \sigma<1$,  and define the algebra $Y$ with this value of $\sigma$, and choice of $V$. 
Assume  $S\in Y$ satisfies, for some $N\ge1$
\EQ{\nn 
\lim_{\eps\to0} \|\eps^{-3}{\chi}(\cdot/\eps)\ast S^{N} - S^{N}\|_{Y} &= 0  \\
\lim_{L\to\infty} \|(1-\hat{\chi}(y/L))S(y)\|_{Y} &=0\nn 
}
Assume  that  $I+\hat{S}(\eta)$ has an inverse in $\B(L^\infty)$ of the form $(I+\hat{S}(\eta))^{-1}= I + U(\eta)$, with  $U(\eta)\in \mc FY$ for all $\eta\in\R^{3}$ uniformly, i.e., 
\EQ{
\nn 
\sup_{\eta\in\R^3} \|U(\eta)\|_{\mc FY}<\infty 
}
Furthermore, let  $\eta\mapsto \hat{S}(\eta)$ be uniformly continuous as a map $\R^{3}\to \B(L^{\infty})$. 
 Then the operator $I+S$ is invertible in the convolution algebra~$Y$. 
\end{theorem}

Lemmas 8.2 and~9.1 in \cite{becsch1} verify that the assumptions of this invertibility theorem hold. Interestingly, Fourier restriction bounds on the resolvent enter  crucially at that stage.  The invertibility pointwise in $\eta$ is precisely the one appearing in the second line of~\eqref{eq:M0}. 
The structure theorems of this section thus depend in an essential way on the newer version of the Agmon-Kato-Kuroda theory based on Fourier restriction (and therefore depend on the non-vanishing curvature of the constant energy surfaces) delineated in Section~\ref{sec:IS}. 

To apply the Wiener theorem above, we return to the key inversion problem~\eqref{keyinv} an conclude that 
\EQ{\label{eq:wichtig}
(I+T_{1+}^{\eps})^{-1}=I-T_{+}^{\eps}, \qquad \sup_{\eps>0} \|T_{+}^{\eps}\|_{Y} <\infty
}
whence 
\EQ{\label{eq:Darstellung}
T_{+}^{\eps} = I - (I+T_{1+}^{\eps})^{-1} & = (I+T_{1+}^{\eps})^{-1} \oast T_{1+}^{\eps} = (I-T_{+}^{\eps})\oast T_{1+}^{\eps} \\
&=  T_{1+}^{\eps} \oast  (I-T_{+}^{\eps})
}
One has the representation formula 
\[
(W_{+}^{\eps}f)(x) = f(x)- \int_{\R^{3}} (\one_{\R^{3}}T_{+}^{\eps})(x,y)f(x-y)\, dy,
\]
where $\frak X_{+}^{\eps}(x,y):= -\one_{\R^{3}}T_{+}^{\eps}$ means the contraction as previously defined. The preceding machinery allows us to conclude that $$\frak X_{+}^{\eps}(x,y)\in X=L^{1}_{y}V^{-1}B_{x}$$
from which the structure formula follows. Indeed, in analogy with the expression~\eqref{eq:key gn} the main term in the definition of the structure function in Theorem~\ref{thm:main1} is 
\EQ{\nn 
h^\eps(x,dy, \omega) = \int_{\R^{3}} g_{1,f^{\eps}_{y'}}^{\eps} (x,d(y-y'),\omega) \, dy'
}
where $f^{\eps}_{y'}(x'):= \frak X_{+}^{\eps}(x',y')$ and $g_{1,f^{\eps}_{y'}}^{\eps}$ refers to the explicit structure function from~\eqref{eq:L}, but with the ``twisted'' potential 
$f^{\eps}_{y'}(x')V(x')$.  The structure function $g$ in Theorems~\ref{thm:main1} and~\ref{thm:main2} is the the sum of $h^{\eps}$ and~$g_{1}^{\eps}$ (for the potential $V$ itself), followed by taking the limit $\eps\to0$.  Note how the Wiener theorem reduces the problem of summing the divergent series of the terms~\eqref{eq:key gn} to the 
inversion problem~\eqref{eq:wichtig} leading to the representation~\eqref{eq:Darstellung}. 

\subsection{  A scaling invariant  condition }

As we already mentioned before, the structure theorems lose a little more than $\frac12$ of a power in terms of decay of~$V$. Ideally, one would wish for a scaling invariant theory. It is perhaps unlikely that this can be achieved in the framework of the spaces $\dot B^{\frac12}$ alone. In~\cite{becsch2} a more complicated  scaling invariant condition on~$V$ was introduced, and a structure theorem for small potentials was obtained in this class. Currently, no analogous version exists for large potentials. 

To briefly describe these results, take a 
Schwartz  potential {   $V$}, and set {    $\trip V\trip:=\|L_V\|_{L^1_{t, \omega}}$}. Recall 
\[
L_V(t,\omega)  =  \int_0^\infty \what{V}(-\tau\omega) e^{\frac{i}{2}t\tau}\, \tau\, d\tau
\]
For any Schwartz function {   $v$} in    $\R^3$ define 
\EQ{\label{eq:Bexotic}
\|v\|_B := \sup_{\Pi} \int_{-\infty}^\infty \trip\delta_{\Pi(t)}\, v(x)\trip \, dt 
}
where {   $\Pi$} is a $2$-dimensional plane through the origin, with all parallel planes {    $\Pi(t)=\Pi+t\vec N$, $\vec N$} being the unit norm to {    $\Pi$}.  Then  with $\psi$ being the usual Littlewood-Paley localizer, one has 
\EQ{\nn 
\|v\|_B\les \sup_{\omega\in \Sph^2} \int_{-\infty}^\infty \sum_{k\in\Z} 2^{\frac{k}{2}} \big\| \psi(2^{-k} x') v(x'+s\omega) \|_{\dot H^{\frac12}(\omega^\perp)}\, ds 
}
The point here is that the right-hand side is formulated in a more accessible way than the implicit norm on the left-hand side. In particular, the right-hand side is finite on Schwartz functions. 

The scaling invariant small potential theorem is the following one:

 \begin{theorem}[\cite{becsch2}]
There exists $c_0>0$ so that  for any real-valued $V$ with 
 $\| V\|_B + \| V\|_{\dot B^{\frac12}}\le c_0$,  there exists  $g(x, y, \omega) \in L^1_\omega \mc M_y L^\infty_x$ 
 with 
 \EQ{\nn 
  \int_{\Sph^2} \|g(x, dy, \omega)\|_{\mc M_y L^{\infty}_x}  \dd \omega \les c_0
} 
such that for any $f \in L^2$ one has the representation formula
\EQ{\nn 
(W_{+} f)(x) &= f(x) + \int_{\Sph^2} \int_{\R^3} g(x, dy, \omega) f(S_\omega x - y)   \dd \omega.
} 
\end{theorem}

In order to prove a large potential analogue, one would need to 
redo all the spectral theory and the Wiener theorem  within the framework of the somewhat exotic $B$-norm from~\eqref{eq:Bexotic}.

\end{document}